\documentclass[10pt]{amsart}

\usepackage{amssymb}
\usepackage{amsmath}
\usepackage{latexsym}

\newtheorem{Theorem}{Theorem}[section]

\newtheorem{Proposition}[Theorem]{Proposition}
\newtheorem{Lemma}[Theorem]{Lemma}
\newtheorem{Corollary}[Theorem]{Corollary}

\newtheorem{Remark}{Remark}[section]

\newcommand{\R}{\ensuremath{\mathbb{R}}}
\newcommand{\N}{\ensuremath{\mathbb{N}}}

\newcommand{\car}{1\hspace{-1.4 mm}1}

\title[Nonlocal HJ equations and dislocations]{NONLOCAL FIRST-ORDER HAMILTON-JACOBI EQUATIONS
MODELLING DISLOCATIONS DYNAMICS}

\author[G. Barles and O. Ley]{Guy Barles \& Olivier Ley\\ \ \\
Laboratoire de Math\'ematiques et  Physique \\ Th\'eorique (UMR CNRS 6083) \\
Facult\'e des Sciences et Techniques \\ Universit\'e de Tours \\
Parc de Grandmont, 37200 Tours, France}

\thanks{The second author is partially supported by a grant ACI JC 1041 
(2002--2005) of the French Ministry of Research.}

\begin{document}

\begin{abstract}
We study nonlocal first-order equations arising in the theory of dislocations. We prove the existence and 
uniqueness of the solutions of these equations in the case of positive and negative velocities, under suitable 
regularity assumptions on the initial data and the velocity. These results are based on new $L^1$-type estimates 
on the viscosity solutions of first-order Hamilton-Jacobi Equations appearing in the so-called ``level-sets approach''. 
Our work is inspired by and simplifies a recent work of Alvarez, Cardaliaguet and Monneau.
\end{abstract}



\maketitle

\section{Introduction}

The starting point of this work and its main motivation is the study of the following type of nonlocal equations 
arising in dislocations' theory \cite{rlf03}
\begin{equation}\label{dyseqn}
u_t=c[\car_{\{u(\cdot ,t)\geq 0\}}]|D u|\quad \hbox {in  }\R^N \times (0,T)\; ,
\end{equation}
where $T > 0$, the solution $u$ is a real-valued function, $u_t$
and $Du$ stand respectively for its time and space derivatives and
$\car_A$ is the indicator function of $A$ for any $A \subset
\R^N.$ For all $\rho \in L^\infty (\R^N)$ or $L^1 (\R^N)$,
$c[\rho]$ is defined by
\begin{eqnarray*}
c[\rho](x,t)=(c_0 * \rho)(x,t)+c_1 (x,t)\ \hbox {in  }\R^N \times (0,T),
\end{eqnarray*}
where $c_0, c_1$ are given functions, satisfying suitable assumptions which are described later on and ``$*$'' 
stands for the usual convolution in $\R^N$ with respect to the space variable $x.$

At first glance, equation (\ref{dyseqn}) looks like equations arising in the so-called ``level-sets approach'' 
to describe the evolution of moving interfaces or domains.  We recall that the level-set
approach was first introduced by Osher and Sethian \cite{os88} for numerical computations and then developed from 
a theoretical point of view by Evans and Spruck \cite{es91} for motion by mean curvature and by Chen, Giga and 
Goto \cite{cgg91} for general normal velocities. We also refer the reader to Barles, Soner and Souganidis \cite{bss93} 
and Souganidis \cite{souganidis97, souganidis95} for different presentations and other results on the level-sets approach.

But, in fact, (\ref{dyseqn}) is not really a level-sets equation, and this for two main reasons : first, in order 
to apply completely viscosity solutions' theory, one would need some monotonicity with respect to the non-local 
dependence in the equation and this would lead here to assume that $c_0 \geq 0$ in $\R^N \times (0,T)$ and  
this assumption is not natural in the dislocations' framework.

Moreover, in the spirit of the level-sets approach, all the level-sets have to be treated in the same way and 
it was remarked by Slep{\v{c}}ev \cite{slepcev03} that, in order to do so, the nonlocal term has to depend typically on sets 
of the form $\{u(\cdot,t)\geq u(x,t)\}$; this is not the case here where the 0-level set plays a particular role.

Finally we point out that the key difficulty in equations like (\ref{dyseqn}) is that, in general, one does not 
expect the indicator function to be continuous as a function of time in $L^1 (\R^N)$ : this is a by-product of 
the well-known ``non-empty interior difficulty'' in the level-sets approach. In particular, to solve (\ref{dyseqn}) 
by approximation turns out either to be very difficult or to lead to very weak formulations. Of course, uniqueness 
is even a more difficult issue and is probably wrong in general.

To the best of our knowledge, the first existence and uniqueness results for (\ref{dyseqn}) in the non-monotone framework 
were obtained by Alvarez, Hoch, Le Bouar and Monneau \cite{ahlm04, ahlm06} : they proved small time results which are
mainly valid for 
graphs but hold without restrictive assumptions on $c_0$ and $c_1$. Then a major breakthrough was made by 
Alvarez, Cardaliaguet and Monneau \cite{acm04} who remarked that, in the situation where $c[\rho]$ is positive for 
any indicator functions (which does not imply that $c_0$ is positive), the existence and uniqueness can be proved 
for any time interval. In order to do so, they use very fine geometrical properties of the moving front 
$\Gamma_t=\{u(\cdot,t)=0\}$ : in particular, they show that, if this front satisfies the interior ball condition at 
time $t=0$, then this property remains true for all time. It gives as a by-product a one-side bound on the curvatures 
of the front and this bound allows to control both the perimeter of $\Gamma_t$ and the volume of  enlarged sets.

Our aim is to simplify the arguments of \cite{acm04} by using a different approach, closer to the spirit of the 
level-set approach: the first step, as in \cite{acm04}, is to obtain fine properties of the solution of the 
standard level-sets equation
\begin{equation}\label{slse}
u_t = c(x,t)|D u|\ \hbox {in  } \R^N \times (0,T) \; ,
\end{equation}
where $c$ is a continuous function, satisfying suitable assumptions and in particular 
$c(x,t) \geq 0$ in $\R^N \times (0,T)$. At this point, it is worth pointing out that we can treat as well the 
case $c(x,t) \leq 0,$ with suitable (and straightforward) adaptations of our arguments and results; we provide 
at the end of Section~\ref{sec:estim} all the needed arguments to do it.

The key result is a $L^1$- estimates on the measure of sets like $\{a \leq
u(\cdot,t)\leq b\}$ where $-\overline\delta \leq a < b \leq \overline\delta$
for some small enough $\overline\delta$. The key difference with \cite{acm04}
is that we use here the classical level-sets approach with
continuous (and even Lipschitz continuous) solutions $u$ while, in \cite{acm04},
just indicator functions are used. In fact, the classical
level-sets solution carries more informations and, roughly
speaking, we replace the fine geometrical estimates of \cite{acm04} on the,
eventually non-smooth, sets $\{u(\cdot,t)=0\}$ by (almost) classical estimates on
$u$ and its derivatives.

To do so, the two key results are the lower bound estimate on $|Du|$ of Ley \cite{ley01} and the more classical 
semiconvex property of the solution of (\ref{slse}). As we mention it above, these two estimates carry the 
necessary informations on the front; maybe we do not obtain as fine estimates as in \cite{acm04} but we obtain 
them in a far simpler way and they are more than enough to study (\ref{dyseqn}), since, in particular they 
imply that $t \mapsto \car_{\{u(\cdot,t)\geq 0\}}$ is
continuous in $L^1(\R^N)$.

Our paper is organized as follows : in Section \ref{sec:estim}, we provide all the necessary results on (\ref{slse}) 
by recalling the classical results. Section \ref{sec:L1} is devoted to the new $L^1$-estimate. 
In Section \ref{sec:ahj}, we describe the 
application to (\ref{dyseqn}) which is obtained by using a classical fixed point arguments for a suitable 
contraction mapping. Several variants exist, either by using a Schauder's fixed point approach and/or an 
approximation argument: each of these approaches has advantages and disadvantages, we made a choice in 
this paper to present one of them, the others will be used in forthcoming works. Finally, in the Appendix, 
we relate in a more precise way the estimates on $u$ we use with some of the geometrical properties obtained 
in \cite{acm04}: this allows the reader to compare more easily the two different approaches and see that they 
are almost equivalent.

\section{Preliminaries on the classical HJ equation: lower-bound gradient estimate,
semiconvexity and front propagation}\label{sec:estim}

We consider the first-order Hamilton-Jacobi equation
\begin{eqnarray}
\left\{
\begin{array}{cc}
u_t = c(x,t)|D u| & {\rm in} \ \R^N\times (0,T), \label{chjb} \\
u(x,0)=u_0(x) & {\rm in} \ \R^N,
\end{array}
\right.
\end{eqnarray}
where $c : \R^N\times [0,T]\to \R$ and $u_0 : \R^N\to \R$ are given continuous
functions, $u_t$ and $D u$ stand, respectively, for the time and space 
derivative of $u$ and $|\cdot |$ is the standard Euclidean norm.

We introduce some assumptions.

\vspace*{2mm}\noindent{\bf (H1)} 
There exist $L_1, L'_1>0$ such that, for all $x,y\in \R^N,$
for all $t\in [0,T],$
\begin{eqnarray}
& & |c(x,t)-c(y,t)|\leq L_1|x-y|, \nonumber \\
& & |c(x,t)|\leq L'_1. \label{cborne}
\end{eqnarray}

\vspace*{2mm}\noindent{\bf (H2)} For all $x\in \R^N, t\in [0,T],$ $c(x,t)\geq 0.$

\vspace*{2mm}
\noindent{\bf (H3)} There exists $\eta_0 >0$ such that,
\begin{eqnarray*}
-|u_0(x)|-|D u_0(x)|+\eta_0 \leq 0 \ {\rm in} \ \R^N \ {\rm in \ the \ viscosity \ sense.}
\end{eqnarray*}

We make some comments about the assumptions.
Note that {\bf (H2)} implies that $p\in\R^N \mapsto c(x,t)|p|$
is convex for every $(x,t)\in\R^N\times [0,T]$ which is a key
assumption here.
When $u_0$ is $C^1,$ {\bf (H3)} expresses that the gradient of
$u_0$ does not vanish  on the set $\{u_0=0\}.$ When $u_0$ is not smooth, 
{\it in the viscosity sense} (see \cite{ley01}) means that $u_0$ is a viscosity subsolution of
equation $-|v(x)|-|D v(x)|+\eta_0\leq 0$ in $\R^N$ or equivalently
that, for all $x\in \R^N$ and $p\in D^+ u_0(x),$
$$
|u_0(x)|+|p|\geq \eta_0,
$$
where $D^+ u_0(x)$ (respectively $D^- u_0(x)$)
denotes the Fr\'echet super-differential (respectively sub-differential) of $u_0$ at $x.$
For viscosity solutions, we refer the reader to \cite{bcd97} and \cite{barles94}.

We say that a function $v:\R^N\to \R$ is  {\it semiconvex} with constant $C>0$ if, for all $x,h\in \R^N,$
\begin{eqnarray} \label{def-sc}
v(x+h)-2v(x)+v(x-h)\geq -C |h|^2.
\end{eqnarray}
We refer to \cite{cs04} for properties of semiconcave and semiconvex functions.
In particular, a semiconvex function is locally Lipschitz continuous and, 
for all $x\in \R^N$ and $p\in D^- v(x),$
\begin{eqnarray} \label{sc-sg}
v(x+h)\geq v(x)+\langle p,h\rangle -\frac{C}{2}|h|^2 
\ \ \ {\rm for \ all} \ h\in \R^N.
\end{eqnarray}
Moreover, a semiconvex function is twice differentiable everywhere
\begin{eqnarray} \label{est-sc}
D^2 v(x)\geq -C Id \ \ \  {\rm for \ a.e.} \ x\in\R^N,
\end{eqnarray}
where $Id$ is the identity matrix in $\R^N.$

\vspace*{2mm}
\noindent{\bf (H4)} $x\in \R^N \mapsto c(x,t)$ is semiconvex with constant $L_2,$
uniformly for all $t\in [0,T].$
\\[2mm]
In the sequel, we denote the essential supremum of $f\in L^\infty (\R^N)$ by
$|f|_\infty.$
\begin{Theorem} \label{thm-borneinf} (\cite{ley01})
\begin{enumerate}
\item[(i)] Under assumption {\bf (H1)}, Equation (\ref{chjb}) 
has a unique continuous viscosity solution $u.$ If $u_0$ is Lipschitz
continuous, then $u$ is Lipschitz continuous and, 
for almost all $x\in \R^N,$ $t\in [0,T],$
$$
|Du(x,t)|\leq {\rm e}^{L_1T}|D u_0|_{\infty}\, , \ \  \ \quad  
|u_t (x,t)|\leq L'_1{\rm e}^{L_1T}|D u_0|_{\infty}\; .
$$
\item[(ii)] Assume that $u_0$ is Lipschitz continuous and that {\bf (H1)}, {\bf (H2)} and {\bf (H3)} hold.
Then there exist $\gamma =\gamma (L_1,L'_1,\eta_0) >0, \eta =\eta (L_1,L'_1,\eta_0) >0$
such that the viscosity solution $u$ of (\ref{chjb})  satisfies in the viscosity sense
\begin{eqnarray}\label{borneinf}
-|u(x,t)|-\frac{e^{\gamma t}}{4}|D u(x,t)|^2+\eta \leq 0 \ { in} \ \R^N\times [0,T] \; .
\end{eqnarray}
\item[(iii)] Assume  that $u_0$ is semiconvex and that {\bf (H1)}, {\bf (H2)}, {\bf (H4)} hold. 
Then $u$ is semiconvex in the $x$-variable
uniformly with respect to $t\in [0,T].$
\end{enumerate}
\end{Theorem}
We refer to \cite{ley01} for proofs of (i)-(ii) and \cite{ley01t} for the proof of (iii). 
We remark that, in (ii), $u$ is Lipschitz continuous
because the assumptions of (i) are satisfied. Therefore $u$ is differentiable a.e. in
$\R^N\times [0,T]$ and (\ref{borneinf}) holds a.e. in $\R^N\times [0,T].$
Part (ii)  gives a lower-bound gradient estimate for $u$ near the front
$\{ (x,t)\in \R^N\times [0,T] : u(x,t)=0\}.$ Indeed, if $|u(x,t)|< \eta /2,$ then 
\begin{eqnarray} \label{bornepresfront}
-|D u(x,t)|\leq -\sqrt{2\eta} e^{-\gamma t/2} <0
\ {\rm in} \ \R^N\times [0,T] 
\end{eqnarray}
in the viscosity sense
hence (\ref{bornepresfront}) holds a.e. in $\R^N\times [0,T].$

We continue by giving an upper-bound for the difference of 
two solutions with different velocity $c.$
\begin{Lemma} \label{traj}
For $i=1,2$, let $u_i \in C(\R^N\times [0,T])$ be a solution of
\begin{eqnarray*}
\left\{
\begin{array}{cc}
(u_i)_t = c_i(x,t)|D u_i | & { in} \ \R^N\times [0,T],\\[2mm]
u_i (x,0)=u_0(x) & { in} \ \R^N,
\end{array}
\right.
\end{eqnarray*}
where $c_i$ satisfies {\bf (H1)} and $u_0$ is Lipschitz continuous. Then, for any $t\in [0,T],$
$$
 |(u_1-u_2)(\cdot ,t)|_\infty
\leq |D u_0|_\infty {\rm e}^{L_1 t}  \int_0^t 
|(c_1-c_2)(\cdot ,s)|_\infty  ds.
$$
\end{Lemma}
\begin{proof}[Proof of Lemma \ref{traj}]
We prefer to focus on the main ideas and so,
part of this proof is formal.
All arguments can be made rigorous using standard properties of viscosity
solutions. Because of Theorem~\ref{thm-borneinf} (i), we have
$$ |Du_i (x,t)|\leq {\rm e}^{L_1T}|D u_0|_{\infty}\quad\hbox{ for  }i=1,2\; ,$$
and therefore in $\R^N\times [0,T]$
$$ (u_1)_t =c_1(x,t)|D u_1 | \leq c_2(x,t)|D u_1 | + |(c_1-c_2)(\cdot ,t)|_\infty{\rm e}^{L_1 t}|D u_0|_\infty \; .$$
It follows that
$$\tilde u_1 (x,t) := u_1 (x,t) - \int_0^t \, |(c_1-c_2)(\cdot ,s)|_\infty{\rm e}^{L_1 s}|D u_0|_\infty ds$$
is a viscosity subsolution of the $u_2$-equation and therefore, by a standard comparison result, $\tilde u_1 \leq u_2$ 
in $\R^N\times [0,T]$, which yields to
$$ u_1 (x,t) - u_2 (x,t) \leq  {\rm e}^{L_1 t}|D u_0|_\infty \int_0^t \, |(c_1-c_2)(\cdot ,s)|_\infty ds\; .$$
The result then follows by exchanging the roles of $u_1$ and $u_2$.
\end{proof}

We turn to an {\it increase principle} for functions satisfying {\bf (H3)}
(see \cite{clsw98} and \cite[Lemma 4.1]{ley01} for similar results). In the sequel,
for any $x\in \R^N$ and $r>0,$ $B(x,r)$ denotes the open Euclidean ball of center $x$
and radius $r$ and $\overline{B}(x,r)$ its closure. 
\begin{Lemma} \label{increasing-princip}
Suppose that $v$ satisfies {\bf (H3)}. Let $\delta <\eta_0/2$ and $x_0\in \{-\delta \leq v\leq \delta\}.$
Then 
\begin{eqnarray} \label{decreas-princ}
\mathop{\rm sup}_{y\in \overline{B}(x_0,2\delta/\eta_0)} v(y)\geq v(x_0)+\delta.
\end{eqnarray}
\end{Lemma}

\begin{proof}[Proof of Lemma \ref{increasing-princip}]
Set $\bar{\eta}=\eta_0/2.$ 
Suppose that (\ref{decreas-princ}) is false. Therefore there exists
$0<\tilde{\delta}<\delta$ such that
\begin{eqnarray} \label{contra-decreas-princ}
\mathop{\rm sup}_{y\in \overline{B}(x_0,\delta/\bar{\eta})} v(y)< v(x_0)+\tilde{\delta}.
\end{eqnarray}
Take $\hat{\delta}, \theta >0$ such that
$0<\tilde{\delta}<\hat{\delta}<\delta$ and $(1+\theta)\hat{\delta} <\delta$ and
set $f(y)=v(y)-\hat{\delta}(\bar{\eta}|y-x_0|/\delta)^{1+\theta}.$
If $y\in \partial B(x_0,\delta/\bar{\eta}),$ using (\ref{contra-decreas-princ}), we have 
\begin{eqnarray*}
f(y)=v(y)-\hat{\delta} < v(x_0)+\tilde{\delta}-\hat{\delta}<f(x_0).
\end{eqnarray*}
Therefore the maximum of $f$ is achieved at $\bar{y}$ lying in the open ball $B(x_0,\delta/\bar{\eta}).$
Moreover, since $-\delta \leq v(x_0)\leq \delta,$ for $\varepsilon>0$ small enough,
$\bar{y}$ belongs to the open set $\{-\delta-\varepsilon <v<\delta +\varepsilon\}$
in which $v$ is a viscosity subsolution of $-|D v|+\bar{\eta}\leq 0$ by {\bf (H3)}.
It follows
\begin{eqnarray*}
\bar{\eta}\leq |D ( \hat{\delta} 
\left( \frac{\bar{\eta} |\cdot -x_0|}{\delta}\right)^{1+\theta})(\bar{y})|
\leq (1+\theta)\frac{\hat{\delta}}{\delta}\bar{\eta} <\bar{\eta}
\end{eqnarray*}
from the choice of $\hat{\delta}$ and $\theta.$
It leads to a contradiction which proves (\ref{decreas-princ}).
\end{proof}

The following lemmas take place in the context of the level-set approach to
front propagation. We refer the reader to \cite{es91}, \cite{cgg91}, \cite{bss93} and \cite{souganidis97,souganidis95} for
details. In few words, in front propagation, we are interested
in the evolution of the set $\Gamma_t= \{u(\cdot, t)=0\}$ which is called the front,
and where $u$ is the continuous viscosity solution of (\ref{chjb}). 
In our case, at least formally, each point $x$ of the front evolves with a normal velocity 
proportional to $c(x,t).$ The level-set approach makes rigorous this evolution
even when $\Gamma_t$ is singular. This approach is based on
the main and surprising result stating that
$\{u(\cdot ,t)= 0\}$ and  $\{u(\cdot ,t)\geq 0\}$ depend only on $\Gamma_0=\{u_0= 0\}$
and $\{u_0\geq 0\}$ (the initial front) and not on the whole function $u_0.$ 
Since (\ref{chjb}) has a ``finite speed of propagation'' property (see \cite[Theorem 6.1]{ley01}),
we have some bounds of the size of the front:
\begin{Lemma} \label{controle-front}
Assume {\bf (H1)}. Suppose $u_0$ is Lipschitz continuous and there exists $R_0>0$
such that $\{ u_0\geq 0\}\subset \overline{B}(0,R_0).$ Let $u$ be the viscosity solution of (\ref{chjb}) 
with initial condition $u_0.$ Then, for all $t\in [0,T],$
\begin{eqnarray*}
\{ u(\cdot ,t)\geq 0 \}\subset \overline{B}(0,R_0+L'_1t),
\end{eqnarray*}
where $L'_1$ is defined in (\ref{cborne}).
\end{Lemma}
\begin{proof}[Proof of Lemma \ref{controle-front}]
The function $u$ is a subsolution of the equation
$$ u_t \leq L'_1|Du|  \quad \hbox{in} \ \R^N\times [0,T], $$
but, for this equation, the Oleinik-Lax formula provides the unique solution and by a standard comparison result, 
we have
$$ u(x,t) \leq \max_{|y-x|\leq L'_1t}\,  u_0 (y) \,  .$$
If $x$ does not belong to $\overline{B}(0,R_0+L'_1t)$, then all point $y$ such that $|y-x|\leq L'_1t$ lies in 
the complementary of the ball $\overline{B}(0,R_0)$ and therefore in the set $\{ u_0< 0\}$. Hence $u(x,t) <0$ and 
the result is proved.
\end{proof}

Moreover, using the lower-bound gradient estimate of Theorem \ref{thm-borneinf}, we
obtain that the front has 0 Lebesgue measure ${\mathcal{L}}^N.$ 
In the sequel, $\car_{A}$ denotes the indicator function of any
measurable set $A.$ 
\begin{Corollary} (\cite[Corollary 5.1]{ley01}) \label{front-mesure0}
Assume {\bf (H1)} and {\bf (H2)}. Suppose that $u_0$ is Lipschitz continuous, 
that {\bf (H3)} holds and that $\{u_0\geq 0\}$ is a compact subset. 
Then, for every $t\in [0,T],$ ${\mathcal{L}}^N (\{u(\cdot ,t)=0\})=0$ and the function 
$t\mapsto \car_{\{u(\cdot ,t)\geq 0\}}$ from $[0,T]$ to $L^1 (\R^N)$ is continuous.
\end{Corollary}
\begin{proof}[Proof of Lemma \ref{front-mesure0}]
As we noticed in (\ref{bornepresfront}), if $|u(x,t)|< \eta /2,$ we have
$-|D u(x,t)|\leq -\sqrt{2\eta} e^{-\gamma t/2} <0$ in the viscosity sense. 
This property is true in $\R^N$ for any fixed 
$t\in (0,T)$ and therefore also almost everywhere in $\R^N$. 
Indeed since the viscosity inequality in $\R^N \times (0,T)$ does not
involve any time-derivative, it is easy to show that it holds in $\R^N$
for any $t$, just by remarking that, for any smooth function $\phi$,
any (strict) local maximum point of $x \mapsto u(x,t) - \phi (x)$ is
approximated by a local maximum point of $(x,s) \mapsto u(x,s) - \phi(x) - (t-s)^2/\varepsilon$ 
where $\varepsilon >0$ is a small parameter devoted to tend to $0$.

Then, for a fixed $t\in (0,T),$ applying the celebrated Stampacchia's result 
(see e.g. \cite[p.84]{eg92}), we know that $D u(\cdot ,t)=0$ almost everywhere on $\{u(\cdot ,t)=0\}$, 
implying that necessarely this set has a zero-Lebesgue measure.

The continuity of the indicator function follows immediately from this property~; let $(x_0,t_0)\in \R^N\times [0,T].$ 
If $u(x_0,t_0)> 0$ (respectively $u(x_0,t_0)< 0$), then, by continuity of $u,$ $u(x_0,t)>0$  (respectively $u(x_0,t)< 0$)
for $t$ close enough to $t_0.$ It follows that $\car_{\{u(\cdot,t)\geq 0\}} (x_0)\to \car_{\{u(\cdot,t_0)\geq 0\}} (x_0)$
as $t\to t_0$ for every $x_0$ such that $u(x_0,t_0)\not= 0.$ But $\{u(\cdot ,t_0)=0\}$ has a zero-Lebesgue measure and 
therefore $\car_{\{u(\cdot,t)\geq 0\}}\to \car_{\{u(\cdot,t_0)\geq 0\}}$ a.e. in $\R^N$ as $t\to t_0.$ And we conclude by the
dominated convergence theorem.
\end{proof}

We conclude this section by mentioning the changes in the above results if, instead of assuming $c  \geq 0$, we 
assume $c \leq 0$. First we point out that Theorem~\ref{thm-borneinf} (i), Lemma~\ref{traj} and \ref{controle-front} 
holds even if $c$ changes sign and therefore these results are independent of the sign of $c$.

Next, if $c \leq 0$, {\bf (H3)} and {\bf (H4}) have to be replaced respectively by\\

\noindent{\bf (H3')} There exists $\eta_0 >0$ such that,
\begin{eqnarray*}
|u_0(x)|+|D u_0(x)|-\eta_0 \geq 0 \ {\rm in} \ \R^N \ {\rm in \ the \ viscosity \ sense.}
\end{eqnarray*}

\noindent{\bf (H4')} $x\in \R^N \mapsto c(x,t)$ is semiconcave with constant $L_2,$
uniformly for all $t\in [0,T].$\\[3mm]
\noindent And under these new assumptions, then, in Theorem~\ref{thm-borneinf}, (\ref{borneinf}) is changed into
\begin{eqnarray*}
|u(x,t)|+\frac{e^{\gamma t}}{4}|D u(x,t)|^2-\eta \leq 0 \ {\rm in} \ \R^N\times [0,T] \; ,
\end{eqnarray*}
while, if $u_0$ is semiconcave, then $u$ is semiconcave with respect to the $x$-variable, uniformly with respect to 
$t \in [0,T]$.

Finally, the {\it increase principle} of Lemma~\ref{increasing-princip} is changed into a {\it decrease principle} which is 
formulated in \cite[Lemma 4.1]{ley01}, while Corollary~\ref{front-mesure0} remains true as a consequence of the new version 
of Theorem~\ref{thm-borneinf} (ii).

In the next section, our estimates rely, roughly and formally speaking, on the fact that the quantity
$$ {\rm div} \left(\frac{c(x,t)D u}{|D u|}\right)\; ,$$
is bounded from above if $(x,t)$ is close enough to the front, i.e. if $u(x,t)$ is small enough. This is based on the 
results (ii) and (iii) of Theorem~\ref{thm-borneinf} in the case $c\geq 0$. It is easy to check that, in the case $c\leq 0$, 
this property is preserved since $u$ is now semiconcave, and the change of sign of $D^2 u$ compensates the change of sign of 
$c$.

\section{Estimates on the measure of small level-sets of the solution
of the HJ equation}
\label{sec:L1}

For every $a<b$ and $\epsilon >0,$ we consider a smooth function $\varphi : \R\to \R^+$ such that
$\varphi =0$ on $(-\infty ,a-\epsilon],$ $\varphi$ is increasing 
on $(a-\epsilon ,a),$
$\varphi =1$ on $[a,b],$ $\varphi$ is decreasing on $(b, b+\epsilon)$
and $\varphi =0$ on $[b+\epsilon, +\infty).$
We choose $\epsilon << b-a$ and $\varphi$ decreasing with respect to $\epsilon$
such that $\varphi \downarrow \car_{[a,b]}$ when $\epsilon\downarrow 0.$
Here $\car_{[a,b]}$ denotes the indicator function of $[a,b].$
Note that we omit to write the dependence of $\varphi$ with respect to 
$a,b$ and $\epsilon$ for the sake of simplicity of notations.

\begin{Proposition} \label{est-int}
Assume {\bf (H1)}, {\bf (H2)}, {\bf (H3)}, {\bf (H4)} and suppose that $u_0$ is 
Lipschitz continuous, semiconvex with constant $L_3$ and
$\{ u_0\geq 0\}$ is a compact subset. Let $-\eta/2 < a-\epsilon <b+\epsilon <\eta /2$ where $\eta$ is 
defined in (\ref{borneinf}) and let $u$ be the continuous viscosity solution of (\ref{chjb}).
Then there exists $L_4=L_4 (L_1, L'_1, L_2, L_3, \eta_0, T)$ such that, for all $t\in [0,T],$
\begin{eqnarray*}
\int_{\R^N} \varphi (u(x,t))dx \leq e^{L_4 t}
\int_{\R^N} \varphi (u_0(x))dx.
\end{eqnarray*}
In particular,
$$ {\mathcal{L}}^N (\{a \leq u(\cdot ,t)\leq b \})\leq e^{L_4 t}{\mathcal{L}}^N (\{a \leq u_0\leq b \})\; .$$
\end{Proposition}

\begin{proof}[Proof of Lemma \ref{est-int}]
The assumptions of Theorem \ref{thm-borneinf} hold. Therefore the solution
$u$ is Lipschitz continuous (with constant $L$) in $\R^N\times [0,T],$ is semiconvex
(with constant $C=C(L_1,L_2,L_3)$) in the $x$-variable and (\ref{chjb}), 
(\ref{est-sc}), (\ref{borneinf}) and (\ref{bornepresfront}) hold a.e. 
in $\R^N\times [0,T].$

To emphasize the main ideas of the proof, we first provide a formal calculation which is justified latter. We have
\begin{eqnarray} \label{derivation1}
\frac{d}{dt}\left( \int_{\R^N} \varphi (u(x,t))dx\right)
= \int_{\R^N} \varphi ' (u(x,t)) u_t (x,t) dx
\end{eqnarray}
for a.e. $t\in [0,T].$
Using Equation (\ref{chjb}), it follows
\begin{eqnarray*}
\int_{\R^N} \varphi ' (u) u_t dx
& = &
\int_{\R^N} \varphi ' (u) c(x,t)|D u|dx \\
&=&
\int_{\R^N} \langle \varphi ' (u) D u , \frac{c(x,t)D u}{|D u|}\rangle dx \\
&=&
\int_{\R^N} \langle D \varphi  (u), \frac{c(x,t)D u}{|D u|}\rangle dx
\end{eqnarray*}
since, from $-\eta/2 < a-\epsilon <b+\epsilon <\eta /2,$ and (\ref{bornepresfront}), we have
$|D u|>\sqrt{2\eta}e^{-\gamma T/2} $ for almost every $(x,t)$ such that 
$\varphi  (u(x,t))\not=0.$
Using an integration by parts, we obtain
\begin{eqnarray}  \label{parpartie1}
\hspace*{-0.5cm} 
\int_{\R^N} \langle D \varphi  (u), \frac{c(x,t)D u}{|D u|}\rangle dx
=
- \int_{\R^N} \varphi  (u) \, {\rm div}( c(x,t) \frac{D u}{|D u|}) dx.
\end{eqnarray}
Applying the lower-bound gradient estimate again and (\ref{est-sc}), we have,
for almost every $(x,t)\in \R^N\times [0,T]$ such that $\varphi  (u(x,t))\not=0,$
\begin{eqnarray} \label{maj-courbure}
\hspace*{-0.8cm} -  {\rm div}( \frac{D u}{|D u|})
= -\frac{1}{|D u|} {\rm trace} \left[ \left( Id-\frac{D u\otimes D u}
{|D u|^2}\right) \nabla^2 u\right]
\leq \frac{e^{\gamma T/2} C}{\sqrt{2\eta}}.
\end{eqnarray}
It gives
\begin{eqnarray} \label{2emediv} 
-  {\rm div}( c(x,t) \frac{D u}{|D u|})
&=& -\langle D c, \frac{D u}{|D u|}\rangle
- c\,  {\rm div}( \frac{D u}{|D u|}) \nonumber \\
&\leq &
L_1 + \frac{e^{\gamma T/2} L'_1 C}{\sqrt{2\eta}}.
\end{eqnarray}
Finally, setting $L_4= L_1+{e^{\gamma T/2} L'_1 C}/{\sqrt{2\eta}},$
we obtain, for a.e. $t\in[0,T],$
\begin{eqnarray*}
\frac{d}{dt}\left( \int_{\R^N} \varphi (u(x,t))dx\right)
\leq L_4 \int_{\R^N} \varphi (u(x,t))dx
\end{eqnarray*}
which yields the result through a classical Gronwall's argument.

It remains to justify (\ref{derivation1}), (\ref{parpartie1}) and the estimates which follow.
From Lemma \ref{controle-front} and since $\{u_0\geq 0\}$ is bounded, $\{u(\cdot ,t)\geq 0\}$ 
belongs to a fixed compact subset for $0\leq t\leq T.$ Moreover, since $u$ satisfies (\ref{borneinf}), 
from Lemma \ref{increasing-princip}, there exists a compact subset $K\in\R^N$
such that, for every $t\in [0,T],$ 
\begin{eqnarray} \label{zonecompacte}
\{ -\eta/2 \leq u(\cdot ,t)\leq \eta/2 \}\subset K.
\end{eqnarray}
Since $-\eta/2 < a-\epsilon <b+\epsilon <\eta /2,$ 
for every $x\in\R^N,$ $0\leq s,t\leq T,$  $s\not= t,$
we then have
\begin{eqnarray*}
\left| \frac{\varphi (u(x,t))-\varphi (u(x,s))}{t-s} \right|
\leq \car_K (x) C_{a,b, \epsilon} L,
\end{eqnarray*}
where $C_{a,b, \epsilon}$ is the Lipschitz constant of $\varphi$
and $L$ is the Lipschitz constant of $u.$
Therefore, we can apply the dominated convergence theorem to obtain
(\ref{derivation1}) when $s\to t.$

The proof of the end of the formal calculation relies on approximation arguments.
We set, for any function $f:\R^N\times [0,T]\to \R,$ $\alpha >0$ 
and $(x,t)\in \R^N\times [0,T],$
\begin{eqnarray*}
f_\alpha (x,t)=(f *\rho_\alpha)(x,t)= \int_{\R^N} f(y,t)\rho_\alpha (x-y)dy
\end{eqnarray*}
where $\rho_\alpha$ is a standard mollifier. Now $u_\alpha$ and $c_\alpha$ are 
$C^\infty$ in space for every $t\in [0,T]$ and $u_\alpha (\cdot,t), c_\alpha(\cdot ,t)
\to u(\cdot ,t),c (\cdot,t)$ as $\alpha\to 0,$ uniformly on compact subsets of $\R^N.$ In particular,
from (\ref{zonecompacte}) and since $-\eta/2 < a-\epsilon <b+\epsilon <\eta /2,$ for $\alpha >0$ small
enough,
\begin{eqnarray*}
\{ a-\epsilon \leq u_\alpha (\cdot ,t)\leq b+\epsilon \}\subset K
\end{eqnarray*}
and therefore $\varphi (u_\alpha)$ has a compact support independent of $\alpha$ and
$t\in [0,T].$
For all $\beta >0,$ we have
$$ \int_{\R^N} \langle D \varphi  (u_\alpha), \frac{c_\alpha (x,t)D u_\alpha}
{\sqrt{|D u_\alpha|^2+\beta}}\rangle dx  =
- \int_{\R^N} \varphi  (u_\alpha) \, {\rm div}( c_\alpha(x,t) \frac{D u_\alpha}
{\sqrt{|D u_\alpha|^2+\beta}}) dx.
$$
From the very definition (\ref{def-sc}), we see that, if $u(\cdot ,t)$ is semiconvex
with constant $C$ then $u_\alpha (\cdot ,t)$ is still semiconvex with the same
constant. Therefore, a similar calculation as (\ref{2emediv}) gives
\begin{eqnarray*}
-  {\rm div}( c_\alpha (x,t) \frac{D u_\alpha }{\sqrt{|D u_\alpha|^2+\beta}})
\leq |D c_\alpha (x,t)| + \frac{C |c_\alpha (x,t)|}
{\sqrt{|D u_\alpha|^2+\beta}}.
\end{eqnarray*}
It follows
$$ \int_{\R^N} \langle D \varphi  (u_\alpha), \frac{c_\alpha (x,t)D u_\alpha}
{\sqrt{|D u_\alpha|^2+\beta}}\rangle dx  \leq 
\int_{\R^N} \varphi  (u_\alpha)
(|D c_\alpha (x,t)| + \frac{C |c_\alpha (x,t)|}
{\sqrt{|D u_\alpha|^2+\beta}})dx.
$$
Now, since $u$ and $c$ are Lipschitz continuous with respect to $x$ (uniformly
with respect to $t\in [0,T]$), $|D u(\cdot,t)|, |D c(\cdot,t)|\in L^1_{\rm loc} (\R^N).$
Thus, $|D u_\alpha (\cdot,t)|\to |D u(\cdot,t)|$ and 
$|D c_\alpha (\cdot,t)|\to |D c(\cdot,t)|$ in $L^1_{\rm loc} (\R^N)$ as $\alpha\to 0$
(see \cite{eg92}).
Sending $\alpha$ to $0,$ we get
\begin{eqnarray*}
\hspace*{-0.5cm} \int_{\R^N} \langle D \varphi  (u), \frac{c (x,t)D u}
{\sqrt{|D u|^2+\beta}}\rangle dx
& \leq & 
\int_{\R^N} \varphi  (u)
(L_1 + \frac{C L'_1}
{\sqrt{|D u|^2+\beta}})dx \\
& \leq & 
\int_{\R^N} \varphi  (u)
(L_1 + \frac{C L'_1}
{\sqrt{2\eta {\rm e}^{-\gamma T} +\beta}})dx
\end{eqnarray*}
and, letting $\beta$ go to $0,$ we conclude as in the formal calculation.

Finally we point out that the second part of the result follows easily by letting $\varepsilon$ tends to $0$.
\end{proof}

\begin{Proposition} \label{est-init}
Under the same assumptions as in Proposition \ref{est-int},
for $-\eta/2 <a<b<\eta/2,$
there exists $L_5=L_5 (L_4, T)$ such that,
for all $t\in [0,T],$
\begin{eqnarray*}
{\mathcal{L}}^N (\{ a \leq u(\cdot ,t) \leq b \}) \leq \frac{L_5 (b-a)}{\eta} 
{\mathcal{L}}^N (B(0,R_0+1))
\end{eqnarray*}
where $R_0$ is such that $\{u_0\geq 0\}\subset B(0,R_0),$
$u$ is the solution of (\ref{chjb}) and $L_4$ is given by Proposition \ref{est-int}.
\end{Proposition}

\begin{proof}[Proof of Lemma \ref{est-init}]
Using Proposition \ref{est-int} and the definition of $\varphi,$ we have
\begin{eqnarray*}
{\mathcal{L}}^N (\{ a \leq u(\cdot ,t) \leq b \})\leq
\int_{\R^N} \varphi (u(x,t))dx \leq e^{L_4 t}
\int_{\R^N} \varphi (u_0(x))dx.
\end{eqnarray*}
Now we estimate the right-hand side of the previous inequality 
proceeding as in Proposition \ref{est-int}.
Since $-\eta/2 <a<b<\eta/2,$
we can take $\epsilon >0$ such that $-\eta/2 <a-\epsilon <b+\epsilon <\eta/2.$
From Lemma \ref{increasing-princip}, we have
\begin{eqnarray*}
\{ a-\epsilon \leq u_0 \leq b+\epsilon \} & \subset & 
\{ u_0 \geq 0 \}+\frac{2 {\rm max}\{0,-a+\epsilon \}}
{\eta}\overline{B}(0,1) \\
& \subset & B(0,R_0+1).
\end{eqnarray*}
Therefore $\varphi (u_0)=0$ outside the ball $B(0,R_0+1).$
The calculation which follows is formal and can be justified as in
Proposition \ref{est-int}. So we skip the complete proof.
We have
\begin{eqnarray*}
\int_{\R^N} \varphi (u_0)dx 
&= & \int_{B(0,R_0+1)} \varphi (u_0)dx \\
&\leq &  \int_{B(0,R_0+1)} \varphi (u_0)\frac{|Du_0|}{\eta}dx \\
&=& \frac{1}{\eta} \int_{B(0,R_0+1)} \langle \varphi (u_0) D u_0,
\frac{D u_0}{|D u_0|}\rangle dx \\
&=& -\frac{1}{\eta} \int_{B(0,R_0+1)} \Phi (u_0)\, 
{\rm div}(\frac{D u_0}{|D u_0|}) dx,
\end{eqnarray*}
where $\Phi$ is the primitive of $\varphi$ which is 0
at $-\infty.$ From the definition of $\varphi,$ we have
$\Phi (r)\leq b-a+2\epsilon$ for all $r\geq a-\epsilon.$
From (\ref{maj-courbure}) (at $t=0$), we get
\begin{eqnarray*}
\int_{\R^N} \varphi (u_0)dx \leq 
\frac{C_0(b-a+2\epsilon)}{\eta} {\mathcal{L}}^N (B(0,R_0+1)),
\end{eqnarray*}
where $C_0$ is the semiconvex constant of $u_0.$
Finally we obtain
\begin{eqnarray*}
{\mathcal{L}}^N (\{ a \leq u(\cdot ,t) \leq b \})\leq
\frac{C_0 e^{L_4 t} (b-a+2\epsilon)}{\eta} {\mathcal{L}}^N (B(0,R_0+1)),
\end{eqnarray*}
which gives the result sending $\epsilon$ to $0.$
\end{proof}

\section{Application to nonlocal HJ modelling dislocation dynamics}
\label{sec:ahj}

In this section, we are going to prove the existence and uniqueness of the solution of the dislocation 
equation (\ref{dyseqn}) by a classical fixed point argument using a suitable contraction map. As we mention 
in the introduction, other types of arguments to prove the existence will be described in a forthcoming paper.

To do so, we consider the Banach space $X = C ([0,T], L^1(\R^N))$ endowed with the norm 
$|\rho|_T = \mathop{\rm sup}_{t\in [0,T]} |\rho (\cdot ,t)|_{L^1}.$
We recall that we are given two continuous functions $c_0,c_1 : \R^N\to\R$ and
define, for any $\rho\in X,$
\begin{eqnarray*}
c[\rho](x,t) =(c_0 * \rho)(x,t)+c_1(x,t)
\end{eqnarray*}
where $\displaystyle (c_0 *\rho)(x,t)=\int_{\R^N} c_0(y-x,t)\rho (y,t) dy$
(note that the convolution is done in space only).
We aim at solving the nonlocal HJ equation
\begin{eqnarray}
\left\{
\begin{array}{cc}
u_t = c[\car_{\{ u(\cdot ,t)\geq 0\} }](x,t)|D u| & {\rm in} \ \R^N\times (0,T), 
\label{nl-hj} \\[2mm]
u(x,0)=u_0(x) & {\rm in} \ \R^N.
\end{array}
\right.
\end{eqnarray}
The main assumption we will use is \\[2mm]
\noindent{\bf (H5)} 
$c_0\in X$ and, for all $x\in\R^N,$ $t\in [0,T],$
$c_1(x,t)\geq |c_0(\cdot,t)|_{L^1}.$

\medskip

We state some regularity properties of $c[\rho].$

\begin{Lemma} \label{op-nl}\ 
\begin{enumerate}
\item[(i)] $c[\rho]$ is well-defined for any $\rho\in X$ and continuous in $\R^N\times [0,T].$ \\
\item[(ii)] Suppose that $c_0, c_1$ satisfy {\rm \bf (H1)}. Then, for any $\rho\in X,$
$x,x'\in\R^N,$ and $t\in [0,T],$
\begin{eqnarray} \label{mauvaise-borne}
& & |c[\rho](x,t)-c[\rho](x',t)|\leq L_1(1+|\rho (\cdot ,t)|_{L^1})|x-x'|, \nonumber \\
& & |c[\rho](x,t)|\leq L'_1(1+|\rho (\cdot ,t)|_{L^1}).
\end{eqnarray}
If $c_0\in X,$ then $|c[\rho](x,t)|\leq |c_0 |_{T} |\rho (\cdot ,t)|_\infty +L'_1.$ \\
In particular, if {\bf (H5)} holds and $|\rho (x ,t)|\leq 1$ for a.e. $(x,t)\in \R^N\times [0,T],$ 
then 
\begin{eqnarray} \label{c0vraimentborne}
|c[\rho](x,t)|\leq |c_0 |_{T}+L'_1.
\end{eqnarray}
\item[(iii)] For any $\rho\in X,$ $x\in\R^N$ and $t\in [0,T],$
$c[\rho](x,t)\geq c_1(x,t)-L'_1 |\rho (\cdot ,t)|_{L^1}.$
If $c_0\in X,$ then $c[\rho](x,t)\geq c_1(x,t)- |c_0(\cdot,t)|_{L^1} |\rho (\cdot ,t)|_\infty.$ \\ 
In particular, if {\bf (H5)} holds and $|\rho (x ,t)|\leq 1$ for a.e. $(x,t)\in \R^N\times [0,T],$ 
then {\bf (H2)} holds for $c[\rho].$ \\
\item[(iv)] If $c_0,c_1$ satisfy {\bf (H4)}, then $c[\rho]$ is semiconvex in $x$ for
any $\rho\in X.$ More precisely, for any $x,h\in \R^N,$ $t\in [0,T],$
\begin{eqnarray*}
\hspace*{-0.5cm}
c[\rho](x-h,t)-2 c[\rho](x,t)+c[\rho](x+h,t)\geq -C(1+|\rho (\cdot ,t)|_{L^1})|h|^2,
\end{eqnarray*}
where $C$ is a semiconvex constant for $c_0$ and $c_1.$
\end{enumerate}
\end{Lemma} 
The proof of this lemma is straightforward so we skip it.

Now, we can state our main result which is equivalent to \cite[Theorem 4.3]{acm04}.
\begin{Theorem} \label{unicite-nl}
Suppose that $c_0,c_1$ satisfy {\bf (H1)}, {\bf (H4)} and {\bf (H5)}.
Assume that $u_0$ is Lipschitz continuous, semiconvex, satisfies {\bf (H3)} and
$\{u_0\geq 0\}$ is a compact subset. 
Then (\ref{nl-hj}) has a unique continuous viscosity solution in $\R^N\times [0,T].$
\end{Theorem}

\begin{proof}[Proof of Lemma \ref{unicite-nl}]
We follow the ideas of the proof of \cite[Theorem 4.3]{acm04} which relies
on a fixed-point theorem. The main difference is that we work with continuous
viscosity solution instead of discontinuous ones.

First notice that, if $c_0,c_1$ satisfy {\bf (H1)}, then
from Lemma \ref{op-nl} and Theorem \ref{thm-borneinf},
for any $\rho\in X,$
\begin{eqnarray}
\left\{
\begin{array}{cc}
u_t = c[\rho](x,t)|D u| & {\rm in} \ \R^N\times (0,T), 
\label{nl-ro} \\[2mm]
u(x,0)=u_0(x) & {\rm in} \ \R^N,
\end{array}
\right.
\end{eqnarray}
has a unique continuous viscosity solution. 
Moreover, we have
\begin{Lemma} \label{cont-L1}
Suppose $c_0,c_1$ satisfy {\bf (H1)} and {\bf (H5)} and
$u_0$ is Lipschitz continuous satisfying {\bf (H3)}.
Let $\rho \in X$ and $u$ be the unique continuous viscosity solution of
(\ref{nl-ro}). Then $t\in [0,T]\mapsto \car_{\{u(\cdot ,t)\geq 0\}}\in L^1(\R^N)$
is continuous.
\end{Lemma}
\begin{proof}[Proof of Lemma \ref{cont-L1}]
Since $c_0,c_1$ satisfy {\bf (H5)}, 
{\bf (H2)} holds for $c[\rho].$ Under {\bf (H1)} and {\bf (H2)}
for $c[\rho]$ and {\bf (H3)} for $u_0$, the results of 
Theorem \ref{thm-borneinf} (i)-(ii), Lemma~\ref{controle-front} and Corollary \ref{front-mesure0} hold
true, providing all the informations needed to prove the result.
\end{proof}

Next we introduce the following set of functions : for $0 \leq \theta \leq \tau \leq T,$ $v\in C(\R^N)$, we denote 
by $X^{\theta, \tau, v}$ the set of functions $\rho\in C([\theta,\tau],L^1(\R^N))$ such that $0 \leq \rho(x,t) \leq 1 $ 
a.e. in $\R^N$ for any $t \in [\theta,\tau]$, $\rho(x,t) = 0$ a.e. for $x \notin \overline{B}(0,R_0+\bar{c}t )$ 
where $\bar{c}$ 
is defined later on, and $\rho (\cdot ,\theta)=\car_{\{v\geq 0\}}.$ This set is clearly a subset of a Banach space 
of the $X$-type which is endowed with the norm
$$ |\rho|_{\theta, \tau} = \mathop{\rm sup}_{t\in [\theta,\tau]} |\rho (\cdot ,t)|_{L^1}.$$

We first define
\begin{eqnarray*}
\begin{array}{cccc}
\Psi : & X^{0, \tau, u_0} & \longrightarrow & X^{0, \tau, u_0} \\
       & \rho & \longmapsto & (t\mapsto \car_{\{u(\cdot ,t)\geq 0\} })_{t \in [0,\tau]},
\end{array}
\end{eqnarray*}
where $u$ is the unique viscosity solution of (\ref{nl-ro}). 

We first show that $\Psi$ is well-defined. The fact that 
$t\mapsto \car_{\{u(\cdot ,t)\geq 0\} }\in C([0,\tau],L^1(\R^N))$ follows directly from 
Lemma~\ref{cont-L1}. Moreover, by assumption, there exists $R_0>0$ such that 
$\{u_0\geq 0\}\subset \overline{B}(0,R_0).$ From {\bf (H5)} and Lemma \ref{op-nl}, for 
all $\rho\in X^{0, \tau, u_0},$ we have the estimate
$|c[\rho]| \leq \bar{c}:=|c_0|_T+L'_1,$ which is independent of $\rho$. By Lemma \ref{controle-front},
\begin{eqnarray*}
\{u(\cdot ,t)\geq 0\}\subset \overline{B}(0,R_0+\bar{c}t )
\ \ \ {\rm for \ all \ } t\in [0,\tau],
\end{eqnarray*}
which is the property required in $X^{0, \tau, u_0}$. It is worth pointing out that this prevents the 
front from blowing-up in finite time (see Remark \ref{exple-explosion} below). In the sequel we denote 
by $M:={\mathcal{L}}^N ( \overline{B}(0,R_0+\bar{c}T) )$.

Next we aim at showing that $\Psi$ is a contraction provided $\tau$ is small enough.

Let $\rho_1,\rho_2\in X^{0, \tau, u_0}$ and let denote by $u_1, u_2$ the viscosity solutions of (\ref{nl-ro})
associated respectively to $\rho_1, \rho_2$ by $\Psi.$ We fix $\overline{\delta} <\eta/2,$ and first choose $\tau$ small
enough in order that
$\displaystyle \mathop{\rm sup}_{t\in [0,\tau]} |(u_1-u_2)(\cdot ,t)|_\infty \leq \overline{\delta}.$
From Lemma \ref{traj}, it suffices to take $\tau$ such that
\begin{eqnarray}\label{choix1}
& & |D u_0|_\infty {\rm e}^{L_1(1+M)T} \int_0^\tau 
|(c[\rho_1]-c[\rho_2])(\cdot ,s)|_\infty ds \nonumber\\
& \leq&  2 |D u_0|_\infty {\rm e}^{L_1(1+M)T} |c_0|_T \tau < \overline{\delta}.
\end{eqnarray}
Set
$\delta = \mathop{\rm sup}_{t\in [0,\tau]} |(u_1-u_2)(\cdot ,t)|_\infty.$
For all $0\leq t\leq \tau,$ we have
\begin{eqnarray*}
\hspace*{-0.5cm} |(\Psi (\rho_1)-\Psi (\rho_2))(\cdot ,t)|_{L^1} &=&
|\car_{\{ u_1(\cdot ,t)\geq 0\}} - \car_{\{ u_2(\cdot ,t)\geq 0\}}|_{L^1} \\
&=&
{\mathcal{L}}^N (\{u_1(\cdot ,t)\geq 0, u_2(\cdot ,t)< 0\}) \\
&& + {\mathcal{L}}^N (\{u_1(\cdot ,t)< 0, u_2(\cdot ,t)\geq 0\}).
\end{eqnarray*}
But, if $x\in  \{u_1(\cdot ,t)\geq 0, u_2(\cdot ,t)< 0\},$ then
$$
-\delta\leq -|(u_1-u_2)(\cdot ,t)|_\infty+u_1(x,t)\leq u_2(x,t)< 0.
$$
Therefore
$$
{\mathcal{L}}^N (\{u_1(\cdot ,t)\geq 0, u_2(\cdot ,t) < 0\})
\leq {\mathcal{L}}^N (\{-\delta \leq u_2(\cdot ,t)< 0\}).
$$
Similarly
$$
{\mathcal{L}}^N (\{u_1(\cdot ,t)< 0, u_2(\cdot ,t)\geq 0\})
\leq {\mathcal{L}}^N (\{-\delta \leq u_1(\cdot ,t) < 0\})
$$
and we obtain, using Proposition \ref{est-init},
\begin{eqnarray*} 
& &  |(\Psi (\rho_1)-\Psi (\rho_2))(\cdot ,t)|_{L^1} \nonumber \\
&\leq& 
{\mathcal{L}}^N (\{-\delta \leq u_1(\cdot ,t)< 0\}) 
 +{\mathcal{L}}^N (\{-\delta \leq u_2(\cdot ,t)< 0\}) \nonumber \\
&\leq& 2 \frac{L_5 \delta}{\eta} {\mathcal{L}}^N (B(0,R_0+1))  \nonumber \\
&\leq& 2 \frac{L_5}{\eta} {\mathcal{L}}^N (B(0,R_0+1))
\mathop{\rm sup}_{t\in [0,\tau]} |(u_1-u_2)(\cdot ,t)|_\infty
\end{eqnarray*}
where $L_5$ is given by Proposition \ref{est-init} replacing $L_1$
by $L_1(1+M).$

We apply Lemma \ref{traj} with $c_i=c[\rho_i],$ $i=1,2.$
It follows
\begin{eqnarray*} 
& |(\Psi (\rho_1)-\Psi (\rho_2))(\cdot ,t)|_{L^1} \nonumber \\
\leq & \displaystyle
\frac{2 L_5 e^{L1(1+M) \tau}}{\eta}
|D u_0|_\infty  {\mathcal{L}}^N (B(0,R_0+1)) \mathop{\rm sup}_{t\in [0,\tau]}
\int_0^t 
|(c[\rho_1]-c[\rho_2])(\cdot ,s)|_\infty ds \nonumber \\
\leq & \displaystyle
C(L_1,L'_1,T,\eta , u_0, M)
\mathop{\rm sup}_{t\in [0,\tau]} \int_0^t|(\rho_1 -\rho_2)(\cdot,s)|_{L^1} \nonumber\\
\leq & \displaystyle C(L_1,L'_1,T,\eta , u_0, M) \tau |\rho_1 -\rho_2|_{0,\tau} .
\end{eqnarray*}
Thus, taking $\tau$ satisfying (\ref{choix1}) and
$\tau\leq (2C(L_1,L'_1,T,\eta , u_0, M))^{-1},$
we obtain 
$|\Psi (\rho_1)-\Psi (\rho_2)|_{0,\tau} \leq 2^{-1} |\rho_1 -\rho_2|_{0,\tau}$
which proves that $\Psi$ is a contraction. Applying the fixed point theorem, we obtain that (\ref{nl-hj})
has a unique solution $\bar{u}$ in $\R^N \times [0,\tau].$ 

Then we repeat the same arguments on the time interval $[\tau ,2\tau]$ by using $X^{\tau ,2\tau, \bar{u}(\cdot , \tau)}$. 
It is worth pointing out that even if, in the above computations, $\tau$ seems to depend on the initial data
(here $\bar{u}(\cdot , \tau)$), all the 
estimates can be shown to be uniform (they depend only on $u_0$ and $R_0$):
indeed a solution on the time interval $[\tau ,2\tau]$ can be seen as a solution 
on $[0 ,2\tau]$ by using the already computed solution on $[0,\tau]$ (it is easy to see that there is no problem for 
$t=\tau$) and therefore all the constants appearing in the upper and lower gradient bounds, the semiconvexity constant 
and the fixed ball $\overline{B}(0,R_0+\bar{c}T )$ which contains the front depend only on the properties of $c_0$, $c_1$ 
and $u_0$ through Lemma~\ref{op-nl} and the results of Section~\ref{sec:estim}.

In order to conclude, we argue by induction, repeating successively the same arguments on $[k\tau ,(k+1)\tau]$, $k\in \N,$ 
till we reach $T$, thus providing a continuous solution on the whole time interval.
\end{proof}

\begin{Remark} \label{exple-explosion}
{\rm 
If we do not have a bound for $|c_0|_{L^1},$ then the only bound for
the speed is (\ref{mauvaise-borne}). This bound is not sufficient to avoid blow-up
phenomenon for the front in finite time. Indeed, consider Equation (\ref{nl-hj})
with $c_0\equiv 1\notin X$ and $c_1\equiv 0.$ This case corresponds
to an evolution of $\Gamma_t:=\{u(\cdot ,t)=0\}$
with a normal velocity proportional to the volume of $\Omega_t:=\{u(\cdot ,t)\geq 0\}.$
Starting with $\Omega_0=B(0,R_0)$ and $\Gamma_0=\partial \Omega_0,$
a straightforward computation shows that $\Omega_t=B(0,R(t))$ and $\Gamma_t=\partial \Omega_t$
where $R(t)$ satisfies the differential equation
\begin{eqnarray*}
\dot R(t)= C_N R(t)^N \ \ \ {\rm for} \ t\geq 0,
\end{eqnarray*}
with $C_N={\mathcal{L}}^N(B(0,1)).$ If $N=1,$ then the evolution is defined for
all $t\geq 0$ by $R(t)=R_0 {\rm e}^{C_Nt}.$ But, when $N\geq 2,$
the evolution is well defined only for $0\leq t <t^*$ by
\begin{eqnarray*}
R(t)=\frac{R_0}{(1-t/t^*)^{1/(N-1)}} \ \ \ {\rm with} \ t^*= \frac{1}{(N-1)C_N R_0^{N-1}}.
\end{eqnarray*}
There is a blow-up at $t^*.$
}
\end{Remark}

\section*{A. \  Semiconvexity, lower-bound gradient estimate
and sets with interior ball condition}

A closed set $E\subset \R^N$ satisfies the
{\it interior ball property} of radius $r>0$ if, at each point $x$
of the boundary $\partial E$ of $E,$ there exists $p\in S^1=\{\xi\in\R^N : |\xi |=1\}$
such that 
\begin{eqnarray*}
\overline{B}(x-rp,r)\subset E.
\end{eqnarray*}
Note that it is equivalent to say that there exists $E_0\subset E$ and $r>0$
such that 
$$
E=E_0+ r \overline{B}(0,1)=\{x+rp : x\in E_0, p\in \overline{B}(0,1)\}.
$$
The link between, on the one hand, the interior ball property and, on the 
another hand, semiconvexity and lower-bound gradient is the following:
\begin{Lemma} \label{equiv-ibp}
Let $E\in \R^N$ be a closed set. Then $E$ satisfies the interior ball property of
radius $r>0$ if and only if there exists a semiconvex function $v:\R^N\to\R$
with semiconvex constant $C$ 
and $\eta_0>0$ such that $v$ satisfies {\bf (H3)} and
\begin{eqnarray} \label{reprE}
{\rm int}(E)=\{ v>0\}, \ \ \ \partial E =\{v=0\}.
\end{eqnarray}
Moreover we have $r\geq \eta_0/C.$
\end{Lemma}
\begin{proof}[Proof of Lemma \ref{equiv-ibp}]
Suppose that $E$ satisfies the interior ball property of radius $r>0.$ Then 
\begin{eqnarray*}
E=E_0+r\overline{B}(0,1) \ \ \ \Leftrightarrow \ \ \
E=\mathop{\bigcup}_{x\in E_0} x+r\overline{B}(0,1).
\end{eqnarray*}
For every $y\in E,$ we define $\phi_y(x)=r^2-|x-y|^2.$
Then $\phi_y$ is semiconvex with constant 2. (Note we can modify $\phi_y$
for $|x-y|\geq 2$ in order to keep $\phi_y$ semiconvex with constant 2 and
to obtain a Lipschitz continuous function with constant 2.)
We set
$$
v(x)= \mathop{\rm sup}_{y\in E_0} \phi_y (x) \leq r^2.
$$
Then $v$ is semiconvex with constant 2 as the supremum of semiconvex functions
with constant 2. Clearly (\ref{reprE}) holds. Moreover, $\phi_y$ is a subsolution
of $-|w|-|D w|+{\rm min}\{2r, r^2\}\leq 0$ in $\R^N.$ Therefore {\bf (H3)}
holds for $v$ with $\eta_0 = {\rm min}\{2r, r^2\} >0$ since $v$ is a
supremum of subsolutions.

We turn to the proof of the inverse implication. Let $x\in \partial E.$
Since $v$ is semiconvex, $v$ is differentiable a.e. in $\R^N$; thus,
there exists a sequence $x_n\to x$ such that $v$ is differentiable at
$x_n$ and, from {\bf (H3)}, ${\rm lim} |D v(x_n)|\geq \eta_0.$
From the upper-semicontinuity of $y\mapsto D^- v(y)$ (see \cite{cs04}), it follows
that there exists $p\in D^- v(x)$ such that $|p|\geq \eta_0.$ Set 
$$
B=B(x+\frac{p\eta_0}{C|p|}, \frac{\eta_0}{C}) \ \ {\rm and} \ \ 
\xi =  y-x-\frac{p\eta_0}{C|p|} \ {\rm for \ all} \ y\in B.
$$
Note that $|\xi|\leq \eta_0/C.$
Using (\ref{sc-sg}) and the fact that $v(x)=0,$ we get
\begin{eqnarray*}
\hspace*{-0.7cm} v(y)=v(x+y-x) & \geq & v(x)+\langle p,y-x\rangle-\frac{C}{2}|y-x|^2 \\
& \geq & \langle p,\xi\rangle +\frac{\eta_0 |p|}{C}-\frac{C}{2}\left( |\xi|^2 +
2\langle \xi ,\frac{\eta_0 p}{C|p|}\rangle+\frac{\eta_0^2}{C^2}\right) \\
&\geq &
\langle p,\xi\rangle\left( 1-\frac{\eta_0}{|p|}\right)
+\frac{\eta_0 |p|}{C} -\frac{\eta_0^2}{2C}-\frac{C|\xi|^2}{2}.
\end{eqnarray*}
Since $|p|\geq \eta_0$ and $|\xi|\leq \eta_0/C,$ it follows
\begin{eqnarray*}
v(y)&\geq&
-|p| |\xi |
\left( 1- \frac{\eta_0}{|p|} \right)
+\frac{\eta_0^2}{C}\left( \frac{|p|}{\eta_0} -1\right)\geq 0
\end{eqnarray*}
which proves that $B\subset E=\{ v\geq 0\}$ and ends the proof.
\end{proof}
\begin{Remark}
{\rm The heuristic idea for the above lemma comes from (\ref{maj-courbure}).
Indeed, if $v$ is sufficiently regular,
$-{\rm div}(D v(x_0)/|D v(x_0)|$ gives the sum of the principal
curvatures of the hypersurface $\{ v= v(x_0)\}.$ }
\end{Remark}

\begin{Remark}
{\rm The kind of equivalence we prove in the lemma was already noticed in earlier 
articles, see for instance Sinestrari \cite{sinestrari04}.
Using this equivalence, Theorem~\ref{thm-borneinf} gives another proof of the preservation of
the interior ball property for attainable sets of control systems (see \cite[Lemma 4.1]{acm04}
for details or Cannarsa and Frankowska \cite{cf05}).} 
\end{Remark}

\section*{B. \ A more precise estimate in terms of perimeter of level-sets}

For any set $E\subset \R^N,$ we define the perimeter ${\rm per}(E)$ by
${\rm per}(E)={\mathcal{H}}^{N-1}(\partial E)$
where ${\mathcal{H}}^{N-1}$ is the $(N-1)$-Hausdorff measure.

Using a result of \cite{acm04} on sets with interior ball property,
we obtain a refinement of Proposition \ref{est-init}.
\begin{Proposition} \label{est-int-per}
Under the assumptions of Proposition \ref{est-int}, suppose that
$-\bar{\eta} < a<b < \bar{\eta}$ where 
$\bar{\eta}= {\rm min}\{ \eta /2, \eta_0/2, \eta_0^2/(4C_0)\}$
with $\eta_0$ given by {\bf (H3)}, $\eta$ given by (\ref{borneinf})
and $C_0$ the semiconvex constant of $u_0.$
Then, for all $t\in [0,T],$
\begin{eqnarray*}
{\mathcal{L}}^N (\{ a \leq u(\cdot ,t) \leq b \}) \leq \frac{3^N e^{L_4 t}(b-a)}{\eta_0}
{\rm per}(\{ u_0\geq b\})
\end{eqnarray*}
where $u$ is the solution of 
(\ref{chjb}) and $L_4$ is given by Proposition \ref{est-int}.
\end{Proposition}

\begin{proof}[Proof of Lemma \ref{est-int-per}]
From Proposition \ref{est-int}, it is enough to find an upper estimate for 
${\mathcal{L}}^N (\{ a\leq u_0 \leq b \}).$
From Lemma \ref{increasing-princip}, for every $x_0\in \{ -\eta_0/2 <a \leq u_0 \leq b< \eta_0/2 \},$
there exists $\bar{y}\in \overline{B}(x_0,2(b-a)/\eta_0)$ 
such that $u_0(\bar{y})\geq b.$ Hence 
$\{ a\leq u_0 \leq b \}\subset \{ b\leq u_0\}+2(b-a) \overline{B}(0,1)/\eta_0.$
From Lemma \ref{equiv-ibp}, the set $\{ b\leq u_0\}$ satisfies the interior ball property of
radius $\eta_0/ (2C_0)$ since $-|D u_0|+\eta_0/2\leq 0$ on $\{u_0 =b\}.$ 
Applying \cite[Lemmas 2.4 and 2.5]{acm04}, we obtain that ${\rm per}(\{ u_0\geq b\}) <+\infty$ and
\begin{eqnarray*}
\hspace*{-0.7cm} {\mathcal{L}}^N (\{ a\leq u_0 \leq b \}) & \leq & 
{\mathcal{L}}^N ( \left( \{ b\leq u_0\}+\frac{2(b-a)}{\eta_0} \overline{B}(0,1) \right)\backslash \{ b\geq u_0\} )\\
& \leq & \frac{\eta_0\, {\rm per}(\{u_0\geq b\})}{2C_0 N}
(\left(1+\frac{4(b-a) C_0}{\eta_0^2}\right)^N-1) \\
& \leq & \frac{2\, 3^{N-1} (b-a)}{\eta_0} {\rm per}(\{u_0\geq b\})
\end{eqnarray*}
since $4(b-a) C_0/\eta_0^2<2$ because of the choice of $\bar{\eta}.$
\end{proof}


\end{document}